\documentclass[letterpaper, 10 pt, conference]{ieeeconf}  

\IEEEoverridecommandlockouts                              
\usepackage{graphicx}
\usepackage{amsmath}
\usepackage{amssymb}
\usepackage{comment}
\usepackage{color}
\usepackage{epstopdf}
\usepackage{float}
\usepackage{mathrsfs}
\usepackage{hyperref}

\usepackage[ansinew]{inputenc}
\usepackage[T1]{fontenc}
\usepackage{amsmath}
\usepackage{amssymb}
\usepackage{comment}
\usepackage{color}
\usepackage{float}
\usepackage{mathrsfs}
\usepackage{hyperref}

\usepackage{amsfonts,latexsym,amsmath,amssymb}
\usepackage[mathscr]{eucal}

\makeatletter
\let\theoremstyle\@undefined                        
\makeatother

\usepackage{amsthm}

\newtheorem{nntheorem}{\bf Theorem}
\newtheorem{nnassumption}{\bf Assumption}
\newtheorem{nndefinition}{\bf Definition}
\newtheorem{nnlemma}{\bf Lemma}
\newtheorem{nncorollary}{\bf Corollary}
\newtheorem{nnproposition}{\bf Proposition}
\newtheorem{nexample}{\bf Example}

\def\sat{\texttt{sat}}

\newenvironment{theorem}
{\begin{nntheorem}\it}
{\end{nntheorem}}
\newenvironment{proposition}
{\begin{nnproposition}\it}
{\end{nnproposition}}
\newenvironment{lemma}
{\begin{nnlemma}\it}
{\end{nnlemma}}

\newenvironment{definition}
{\begin{nndefinition}\it}
{\end{nndefinition}}
\newenvironment{assumption}
{\begin{nnassumption}\it}
{\end{nnassumption}}
\newenvironment{example}
{\begin{nnexample}\rm}
{\end{nnexample}}
\newenvironment{nnexample}
{\begin{nexample}\rm}{\end{nexample}}

\newtheorem{nnremark}{\bf Remark}
\newenvironment{remark}{\begin{nnremark}}{\hfill \hspace*{1pt}\hfill $\circ$\end{nnremark}}

\newenvironment{proofof}{{\bf Proof of}}{\hfill \hspace*{1pt}\hfill $\Box$}

\newcommand{\stopSwann}{\color{black}}

\def\RR{\mathbb{R}}

\def\epsilon{\varepsilon}

\title{\LARGE \bf
On ISS-Lyapunov functions for infinite-dimensional linear control systems subject to saturations}

\author{Swann Marx$^{1}$, Yacine Chitour$^{2}$ and Christophe Prieur$^{3}$
\thanks{$^{1}$Swann Marx is with LAAS-CNRS, Universit\'e de Toulouse, CNRS, 7 avenue du colonel Roche, 31400, Toulouse, France
        {\tt\small marx.swann@gmail.com}.}%
\thanks{$^{2}$Yacine Chitour is with  Laboratoire des Signaux et Syst\`emes
(L2S), CNRS - CentraleSupelec - Universit\'e Paris-Sud, 3, rue Joliot Curie, ´
91192, Gif-sur-Yvette, France,
{\tt\small yacine.chitour@lss.supelec.fr}.}%
\thanks{$^{3}$Christophe Prieur is with Univ. Grenoble Alpes, CNRS, Gipsa-lab, F-38000 Grenoble, France,
   {\tt\small christophe.prieur@gipsa-lab.fr}.}
  \thanks{This research was partially supported by the iCODE Institute, research project of the IDEX Paris-Saclay, by the Hadamard Mathematics LabEx (LMH) through the grant number ANR-11-LABX-0056-LMH in the "Programme des Investissements d'Avenir" and by the European Research Council (ERC)
through an ERC-Advanced Grant for the TAMING project (grant agreement 66698).}
}

\begin{document}

\maketitle

\begin{abstract}
This article deals with the derivation of ISS-Lyapunov functions for infinite-dimensional linear systems subject to saturations. Two cases are considered: 1) the saturation acts in the same space as the control space; 2) the saturation acts in another space, especially a Banach space. For the first case, an explicit ISS-Lyapunov function can be derived. For the second case, we can only ensure the existence of an ISS-Lyapunov function.  
\end{abstract}



\section{Introduction}\label{introduction}

In recent decades, a great deal of effort has been dedicated to the development of tools for the analysis of systems with a saturated control in a finite-dimensional framework (see \cite{tarbouriech2011book_saturating}, for a nice introduction to the topic as well as the main techniques used). It is well known, in the finite-dimensional framework, that without any Lyapunov stability property for the open-loop system, the origin cannot be globally asymptotically stable when closing the loop with a saturated feedback-law. Moreover, even if the open-loop system is Lyapunov stable, simply saturating a stabilizing feedback can lead to undesirable behavior for the asymptotic stability of the closed-loop system (see e.g., \cite{Fu}). In these cases, the feedback-law has to be designed taking into account the nonlinearity (see e.g. the famous nested saturation solution given in \cite{teel1992globalsaturation} for the chain of integrators and generalized in \cite{sussmann1991saturation}, and \cite{met96} for another solution based on optimization). 

The second aim of the studies of saturated control system deals with the robustness issue: once a saturated feedback law is designed to stabilize a control system, one generally considers the closed-loop system subject to disturbances. It turns out that the feedback law provided by minimizing time (with bounded inputs) and the nested saturated feedback law are not robust (as defined previously) (see e.g., \cite{rao:mag:2001naive}). However \cite{Sab0} provides a robust feedback law, which is a modification of the feedback law furnished in \cite{met96}.  Other easier implementable robust feedback laws have been proposed, at least for the chain of integrators (see e.g., \cite{CHL15}). 

At the best of our knowledge, analysis of infinite-dimensional systems subject to saturations started with \cite{slemrod1989mcss} and \cite{seidman2001note} and it is only very recently that appeared other works such as  \cite{mcpa2015kdv_saturating} and \cite{prieur2016wavecone}. In  \cite{slemrod1989mcss} and \cite{seidman2001note}, global asymptotic stability of the closed-loop system is tackled using nonlinear semigroup theory. Moreover \cite{mcpa2015kdv_saturating} and \cite{prieur2016wavecone} use similar results, but for a wider class of saturations. In \cite{mcpa2017siam}, a nonlinear partial differential equation, namely the Korteweg-de Vries equation, is considered and some Lyapunov arguments are used to conclude on the global asymptotic stability of the closed-loop system. Finally, some nonlinear open-loop abstract control systems are studied in \cite{map2017mcss}, where a LaSalle Invariance Principle is applied to prove global asymptotic stability. 

A system subject to disturbances is said to be input-to-state stable (for short ISS) if its state is bounded by the norm of the disturbances. It implies that the state is bounded if the disturbances are bounded, and converges to $0$ if the disturbances goes to $0$. Such a notion, introduced in \cite{sontag1989smooth} for finite-dimensional systems, is an important tool to design robust feedback laws or even observers. This kind of properties has been studied for finite-dimensional control systems subject to saturations (see e.g., \cite{azouit2016strong}). Note that the ISS notion has been extended for infinite-dimensional systems (see e.g., \cite{prieur2011iss}, \cite{karafyllis2016iss} and \cite{MiW17b}).

The aim of this paper is to study linear infinite-dimensional systems subject to saturations and disturbances, therefore to find ISS-Lyapunov functions for such systems. Two cases are considered: 1)  the saturation acts in the same space as the control space; 2) the saturation acts in another space, especially a Banach space. For the first case, we provide an ISS-Lyapunov function inspired by \cite{liu1996finite}. For the second case, we provide the ISS property, without giving any explicit Lyapunov function.

The paper is organized as follows. In Section \ref{sec_strict_lyap}, the problem is stated, some useful definitions are introduced and the main result is provided. In Section \ref{sec_proof}, some definitions and technical lemmas are given together with the proof of the main result. Some numerical simulation are collected in Section \ref{sec:ex} for the linearized Korteweg-de Vries equation. Finally, Section \ref{sec_conclusion} collects some concluding remarks.

\section{Problem statement and main results} 
\label{sec_strict_lyap}

\subsection{Some remarks on the linear system}

Let $H$ and $U$ be real Hilbert spaces equipped with the scalar product $\langle \cdot,\cdot\rangle_H$ and $\langle \cdot,\cdot\rangle_U$, respectively. Let $A:D(A)\subset H\rightarrow H$ be a (possibly unbounded) linear operator whose domain $D(A)$ is dense in $H$. Suppose that $A$ generates a strongly continuous semigroup of \emph{linear} contractions denoted by $e^{tA}$. Let $B: U\rightarrow H$ be a bounded operator. 
The operators $A^\star$ and $B^\star$ denote the adjoint operators of $A$ and $B$, respectively.
The aim of this paper is to study the following abstract control system
\begin{equation}
\label{abstract-system}
\left\{
\begin{split}
&\frac{d}{dt}z = Az+Bu,\\
&z(0)=z_0.
\end{split}
\right.
\end{equation}
Suppose that we have the following output
\begin{equation}
y=B^\star z.
\end{equation}
Therefore, in this case, a natural feedback law for \eqref{abstract-system} is $u=-B^\star z$. One needs the following assumption. 
\begin{assumption}[Linear feedback law]
\label{linearfeedbacklaw}

The feedback law $u=-B^\star z$ makes the origin of the closed-loop system
\begin{equation}
\label{lclosed-loop}
\left\{
\begin{split}
&\frac{d}{dt}z = (A-BB^\star)z,\\
&z(0)=z_0.
\end{split}
\right.
\end{equation}
globally exponentially stable. 
\end{assumption}
It is easy to see that, under Assumption \ref{linearfeedbacklaw}, $A-BB^\star$ generates a strongly continuous semigroup of contractions (see \cite{seidman2001note} for a precise proof of this result).  In the remaining sections of this paper, the following notation will be used.
\begin{equation}\label{eq:def:tildeA}
\tilde{A} := A-BB^\star
\end{equation}
Let us denote the strongly continuous semigroup of contractions generated by $\tilde{A}$ by $e^{t\tilde A}$.

From \cite[Theorem 5.1.3, page 217]{curtain1995introduction}, Assumption \ref{linearfeedbacklaw} holds if and only if there exist a self-adjoint and definite positive operator $P\in\mathcal{L}(H)$ and a positive value $C$ such that
\begin{equation}
\label{lyapunov-linear-equation}
\langle \tilde A z,Pz,\rangle_H + \langle P z,\tilde Az\rangle_H\leq-C\Vert z\Vert^2_H,\quad \forall z\in D(\tilde{A}),
\end{equation}
where $D(\tilde{A})$ is the domain of $\tilde A$ (by (\ref{eq:def:tildeA}) and the boundedness of $B$, it coincides with the domain of $A$). 
In other words, this inequality implies that there exists a Lyapunov function for the closed-loop system \eqref{lclosed-loop}, given by
\begin{equation}
\label{lyapunov-linear}
V(z):=\langle Pz,z\rangle_H.
\end{equation}

\subsection{Saturation functions}
In this section, the case where the feedback law is bounded is studied. As explained in introduction, the feedback laws rely on the use of a saturation function defined next.

\begin{definition}[Saturation functions on $S$]
\label{def-sat}
Let $S$ be a real Banach space equipped with the norm $\Vert \cdot \Vert_S$. Assume moreover that $(U,S)$ is a rigged Hilbert space. In other words, $S$ is a dense subspace of $U$ and the following inclusions hold
\begin{equation}
S\subseteq U\subseteq S^\prime.
\end{equation}
In particular, the duality pairing between $S$ and $S^\prime$ is compatible with the inner product on $U$, in the sense that
\begin{equation}
(v,u)_{S^\prime\times S}=\langle v,u\rangle_U,\quad \forall u\in S\subset U,\quad \forall v\in U=U^\prime\subset S^\prime.
\end{equation} A function $\sigma:S\rightarrow S$ is said to be a \emph{saturation function on $S$ admissible for $U$} if it satisfies the following properties
\begin{itemize}
\item[1.] \label{bounded} For any $s\in S$, the following holds
\begin{equation}
\Vert \sigma(s)\Vert_S\leq 1;
\end{equation}
\item[2.] For any $s,\tilde{s}\in S$, the following holds
\begin{equation}
\langle \sigma(s)-\sigma(\tilde{s}),s-\tilde{s}\rangle_U\geq 0;
\end{equation}
\item[3.] For any $s,\tilde{s}\in S$, there exists a positive value $k$ such that the following holds
\begin{equation}
\Vert \sigma(s)-\sigma(\tilde{s})\Vert_U\leq k\Vert s-\tilde{s}\Vert_U; 
\end{equation}
\item[4.] For any $s\in U$, the following holds
\begin{equation}
\Vert \sigma(s) - s\Vert_{S^\prime}\leq \langle \sigma(s),s\rangle_U.
\end{equation}
\item[5.] For any $s,\tilde s \in U$, there exists $C_0$ such that the following holds
\begin{equation}
\langle s,\sigma(s+\tilde s)-\sigma(s)\rangle_U \leq C_0\Vert \tilde s\Vert_U.
\end{equation}
\end{itemize}
\end{definition}

\begin{example}[Examples of saturations]
Here are examples of saturations borrowed from \cite{tarbouriech2011book_saturating} or \cite{slemrod1989mcss}.
\begin{itemize}
\item[1.] Suppose that $S:=L^\infty(0,1)$ and $U=L^2(0,1)$. In this case, $S^\prime$ is the space of Radon measures. It contains strictly the space $L^1(0,L)$. In particular, the norm of $S^\prime$ is identified by the norm of $L^1(0,L)$. Hence, $(U,S)$ is a rigged Hilbert space. 

The classical saturation, which is used in practice, can be defined as follows:
\begin{equation}
\label{satLinfty}
\begin{split}
L^\infty(0,1) &\rightarrow L^\infty(0,1),\\
s &\mapsto \mathfrak{sat}_{L^\infty(0,1)}(s),
\end{split}
\end{equation}
where 
\begin{align}
&\mathfrak{sat}_{L^\infty(0,1)}(s)(x) = \sat_{\RR}(s(x))
\end{align}
and
\begin{align}
\sat_{\RR}(s):=\left\{
\begin{array}{rl}
-1 &\text{ if } s\leq -1,\\
s &\text{ if } -1\leq s\leq 1,\\
1 &\text{ if } s\geq 1,
\end{array}
\right.
\end{align} 
All the items of Definition \ref{def-sat} are proved in \cite{bible_khalil} and \cite{tarbouriech2011book_saturating}.
\item[2.] Suppose that $S:=U$. The saturation, studied in \cite{slemrod1989mcss}, \cite{lasiecka2002saturation} and \cite{mcpa2017siam}, is defined as follows, for all $s\in U$
\begin{equation}
\label{sath}
\mathfrak{sat}_U(s):=\left\{
\begin{array}{rl}
&\hspace{-1.1cm}s\hspace{0.9cm}\text{ if }\Vert s\Vert_{U}\leq 1,\\
\frac{s}{\Vert s\Vert_{U}}&\text{ if } \Vert s\Vert_{U}\geq 1.
\end{array}
\right.
\end{equation}
This is obvious that this operator satisfies Item 1. of Definition \ref{def-sat}. In \cite{seidman2001note}, this operator is proved to be $m$-dissipative, which implies that it satisfies Item 2. of Definition \ref{def-sat}. The fact that this operator is globally Lipschitz is proven in \cite{slemrod1989mcss}. Hence, this operator satisfies Item 3. 
\end{itemize}
Let us prove that $\mathfrak{sat}_U$ satisfies Item 4. of Definition \ref{def-sat}.

\textbf{Case $1$: $\Vert s\Vert_U \leq 1$.}

In this case, one has 
$\mathfrak{sat}_U(s)=s$ and 
\begin{equation}
\Vert \mathfrak{sat}_U(s)-s\Vert_U=0\leq \langle \mathfrak{sat}_U(s),s\rangle_U =\Vert s\Vert^2_U.
\end{equation}

\textbf{Case $2$: $\Vert s\Vert_U\geq 1$.}

In this case, the following holds
\begin{equation}
\mathfrak{sat}_U(s) = \frac{s}{\Vert s\Vert_U}
\end{equation}
Hence, it is immediate  to see that, for all $s\in U$ satisfying $\Vert s\Vert_U\geq 1$, one has 
\begin{align}
\Vert \mathfrak{sat}_U(s)-s\Vert_U \leq \langle \mathfrak{sat}_U(s),s\rangle_U\leq \Vert s\Vert_U^2.
\end{align}
Hence, $s\in U\rightarrow \mathfrak{sat}(s)\in U$ satisfies Item 4 of Definition \ref{def-sat}.

Let us prove that $s\in U\rightarrow \mathfrak{sat}(s)\in U$ satisfies Item 5. of Definition \ref{def-sat}. The proof is inspired by \cite{slemrod1989mcss}.

\textbf{Case $1$: $\Vert s\Vert_U\leq 1$, $\Vert s+\tilde s\Vert_U\leq 1$.} 

In this case, the following holds
\begin{align*}
\langle s,\mathfrak{sat}_U(s+\tilde{s})-\mathfrak{sat}_U(s)\rangle_U&=\langle s,\tilde{s}\rangle_U\\
&\leq \Vert \tilde{s}\Vert_U
\end{align*}
\textbf{Case $2$: $\Vert s\Vert_U\geq 1$, $\Vert s+\tilde s\Vert_U\geq 1$}
\begin{align*}
\langle s,\mathfrak{sat}_U(s+\tilde{s})-\mathfrak{sat}_U(s)\rangle_U=&\left\langle s,\frac{s+\tilde{s}}{\Vert s+\tilde{s}\Vert_U}-\frac{s}{\Vert s\Vert_U}\right\rangle_U\\
\leq &\Vert s\Vert_U\left\Vert\frac{s+\tilde{s}}{\Vert s+\tilde{s}\Vert_U}-\frac{s}{\Vert s\Vert_U} \right\Vert_U
\end{align*}
We have
\begin{align*}
\left\Vert\frac{s+\tilde{s}}{\Vert s+\tilde{s}\Vert_U}-\frac{s}{\Vert s\Vert_U} \right\Vert=&\left\Vert\frac{(s+\tilde{s})\Vert s\Vert_U -  s\Vert s+\tilde{s}\Vert_U}{\Vert s+\tilde{s}\Vert_U\Vert s\Vert_U}\right\Vert_U\\
\leq &\frac{1}{\Vert s\Vert_U\Vert s+\tilde{s}\Vert_U}2\Vert s+\tilde{s}\Vert_U\Vert s+\tilde{s}-s\Vert_U\\
\leq &  \frac{2}{\Vert s\Vert_U}\Vert \tilde{s}\Vert_U
\end{align*}
Therefore, one has
\begin{equation}
\langle s,\mathfrak{sat}_U(s+\tilde{s})-\mathfrak{sat}_U(s)\rangle_U\leq 2\Vert \tilde{s}\Vert_U
\end{equation}
\textbf{Case 3: $\Vert s\Vert\leq 1$, $\Vert s+\tilde{s}\Vert_U\geq 1$}

Using the fact that $s\mapsto \mathfrak{sat}_U(s)$ is globally Lipschitz with constant $3$ yields
\begin{align*}
\langle s,\mathfrak{sat}_U(s+\tilde{s})-\mathfrak{sat}_U(s)\rangle_U\leq &\Vert s\Vert_U\Vert \mathfrak{sat}_U(s+\tilde{s})-\mathfrak{sat}_U(s)\Vert_U\\
\leq &3\Vert \tilde{s}\Vert_U.
\end{align*}
\textbf{Case 4: $\Vert s\Vert\geq 1$, $\Vert s+\tilde{s}\Vert_U\leq 1$}
In this case, one has
\begin{align*}
\langle s,\mathfrak{sat}_U(s+\tilde{s})-\mathfrak{sat}_U(s)\rangle_U = &\left\langle s,s+\tilde{s}-\frac{s}{\Vert s\Vert_U} \right\rangle\\
\leq &\Vert s\Vert_U\left\Vert \frac{(s+\tilde{s})\Vert s\Vert-s}{\Vert s\Vert_U}\right\Vert_U\\
\leq & \Vert s+\tilde{s}\Vert(3\Vert s+\tilde{s}-s\Vert_U)\\
\leq & 3\Vert \tilde{s}\Vert_U
\end{align*}
Therefore $\mathfrak{sat}_U$ satisfies all items of Definition \ref{def-sat} for an admissible saturation map.
\end{example}

\subsection{Discussion on the well-posedness}

\subsubsection{Case without disturbance}

We aim at studying the following closed-loop system
\begin{equation}
\label{nlclosed-loop}
\left\{
\begin{split}
&\frac{d}{dt}z = Az-B\sigma(B^\star z):=A_\sigma z,\\
&z(0) = z_0,
\end{split}
\right.
\end{equation}
where the operator $A_\sigma$ is defined as follows
\begin{equation}
\begin{split}
A_\sigma: D(A_\sigma)\subset H&\rightarrow H\\
z&\mapsto A_\sigma z,
\end{split}
\end{equation}
with $D(A_\sigma)=D(A)$. 

From \cite[Lemma 2.1., Part IV, Page 165.]{showalter1997monobanach}, we have the following proposition
\begin{proposition}[Well-posedness]\label{theo:well:posed}
The operator $A_\sigma$ generates a strongly continuous semigroup of contractions, which we denote by $(W_\sigma(t))_{t\geq 0}$. Hence, there exists a unique strong solution to \eqref{nlclosed-loop} and, for all $z_0\in D(A_\sigma)$, the following functions
\begin{equation}
\label{eq:decrease}
t\mapsto \Vert W_\sigma(t)z_0\Vert_H,\quad t\mapsto \Vert A_\sigma W_\sigma(t)z_0\Vert_H
\end{equation}
are nonincreasing. 
\end{proposition}

As far as the well-posedness of the closed-loop system \eqref{nlclosed-loop} is concerned, the interested reader can read \cite{seidman2001note}, \cite{slemrod1989mcss} or \cite{map2017mcss}. The case under consideration in Proposition \ref{theo:well:posed}, the proof is easier than in the previous work, since the nonlinear operators $\sigma$ considered in Definition \ref{def-sat} are globally Lipschitz and monotone. The proof of Proposition \ref{theo:well:posed} is omitted to focus, in this conference paper, on the main result, namely the stability and ISS properties of a nonlinear infinite-dimensional systems, and its proofs.

\subsubsection{Case with disturbance}

In this paper, we will also consider the case with a perturbation. The system under consideration is the following
\begin{equation}
\label{sys-dis}
\left\{
\begin{split}
&\frac{d}{dt}z = Az-B\sigma(B^\star z +d),\\
&z(0)=z_0,
\end{split}
\right.
\end{equation}
where $d:\mathbb{R}_{\geq 0}\rightarrow U$ is the perturbation. Note that $d$ belongs to $\mathfrak{U}:=L^2((0,\infty);U)$. 

For this system, we study only weak solution, that can be written as follows, for every $z_0\in H$.
\begin{equation}
z(t) = e^{tA}z_0 - \int_0^t e^{(t-s)A}B\sigma(B^\star z(s) + d(s))ds.
\end{equation}

Let us introduce two notions related to the ISS property.
\begin{definition}[ISS property]\label{def:iss}
\item[1. ] The origin of \eqref{sys-dis} is said to be ISS with respect to $d$ if there exist a class $\mathcal{KL}$ function $\alpha$ and a class $\mathcal{K}_\infty$ function $\rho$ such that, for any weak solution to \eqref{sys-dis}
\begin{equation*}
\Vert z\Vert_H\leq \alpha(\Vert z_0\Vert_H,t) + \rho(\Vert d\Vert_{\mathfrak{U}}),\quad \forall t\geq 0.
\end{equation*}
\item[2.] A positive definite function $V:H\rightarrow \mathbb{R}_{\geq 0}$ is said to be an ISS-Lyapunov function with respect to $d$ if there exist two class $\mathcal{K}_\infty$ functions $\alpha$ and $\rho$ such that, for any weak solution to \eqref{sys-dis}
\begin{equation*}
\frac{d}{dt}V(z)\leq -\alpha(\Vert z\Vert) + \rho(\Vert d\Vert_U).
\end{equation*}
\end{definition}
\begin{remark}
\begin{itemize}
\item Note that the Lyapunov function $V$ might be non-coercive. In this case, in general, Item 2. of Definition \ref{def:iss} does not imply the ISS property (see e.g., \cite{MiW17b}). Another property is needed, which is the following: for any positive values $C$ and $\tau$, it holds that
\begin{equation}
\label{BRS}
\sup\lbrace \Vert z\Vert\mid \Vert z_0\Vert_H\leq C,\: \Vert d\Vert_H\leq C,\: t\in [0,\tau]\rbrace<\infty.
\end{equation}
In the special case of system \eqref{sys-dis}, this property holds. Checking it reduces to performing the derivative of the function $\Vert z\Vert^2_Z$ along the solutions to \eqref{sys-dis} and to using Item 5. of Definition \ref{def-sat}. Indeed, one obtains
\begin{equation}
\frac{d}{dt}\Vert z\Vert_H^2\leq C_0\Vert d\Vert_U,
\end{equation}
which implies
\begin{equation}
\Vert z\Vert^2_H\leq \Vert z_0\Vert^2_H + C_0\Vert d\Vert_{\mathfrak{U}}.
\end{equation}
This latter inequality implies the property defined by \eqref{BRS}. 
\item It is not sure yet that Item 1. implies the existence of an ISS Lyapunov function (see e.g, \cite{MiW16b} and \cite{MiW17b}). 
\end{itemize}
\end{remark}

\subsection{Input-to-state stability result}

We are now in position to state our main result. Here is its statement.
\begin{theorem}[ISS result]

\label{thm-abstract-lyap}
Suppose that Assumption \ref{linearfeedbacklaw} holds and let $P\in\mathcal{L}(H)$ be the self-adjoint and positive operator $P\in\mathcal{L}(H)$ satisfying \eqref{lyapunov-linear-equation}. Let $\sigma$ be a saturation function satisfying all the items of Definition \ref{def-sat}.
Then, the following holds. 
\begin{itemize}
\item[1.] If $S=U$, there exists an ISS-Lyapunov function for \eqref{sys-dis}.
\item[2.] If $S\neq U$, assume that the following holds, for any $s\in D(A)$,
\begin{equation}
\label{continuous-embedding}
\Vert B^\star s\Vert_S \leq c_S \Vert s\Vert_{D(A)}.
\end{equation}
Then, the origin of \eqref{sys-dis} is ISS with respect to $d$.
\end{itemize} 
\end{theorem}

\section{Proof of the main result}

\label{sec_proof}

\subsection{Definitions and techninal lemmas}

Before proving Theorem \ref{thm-abstract-lyap}, let us give some useful definitions and results to prove Theorem \ref{thm-abstract-lyap}, especially the case where $S\neq U$. 

The following result links the disturbed system \eqref{sys-dis} to the undisturbed one \eqref{nlclosed-loop}. 

\begin{lemma}
\label{lemma-estimation}
If $0$ is globally asymptotically stable for the following abstract control system, for any $z_0\in H$
\begin{equation}
\label{undisturbed}
\left\{
\begin{split}
&\frac{d}{dt}z = Az-B\sigma(B^\star z),\\
&z(0)=z_0,
\end{split}
\right.
\end{equation}
then $0$ is ISS with respect to $d$ for the following abstract control system, for any $z_0\in H$
\begin{equation}
\label{disturbed}
\left\{
\begin{split}
&\frac{d}{dt}z^d = Az^d-B\sigma(B^\star z^d+d),\\
&z^d(0)=z_0.
\end{split}
\right.
\end{equation}
\end{lemma}


Let us prove this lemma.

\begin{proofof}\textbf{ Lemma \ref{lemma-estimation}:}
Let us introduce $\tilde{z}:=z^d-z$. It satisfies the following abstract system
\begin{equation}
\label{abstract-difference}
\left\{
\begin{split}
&\frac{d}{dt}\tilde z = A\tilde z+B\sigma(B^\star z)-B\sigma(B^\star z^d + d),\\
&z(0)=0.
\end{split}
\right.
\end{equation}
Moreover, we compute 
\begin{align*}
\frac{d}{dt}\Vert \tilde z \Vert_H^2 = &\langle A\tilde z,\tilde z\rangle_H + \langle \tilde z ,A\tilde z\rangle_H\\
&+2 \langle B(\sigma(B^\star z)-\sigma(B^\star z^d + d)),\tilde z\rangle_H\\
\leq & \Vert \sigma(B^\star z)-\sigma(B^\star z^d +d)\Vert_U\Vert B^\star \tilde z\Vert_U\\
\leq & k\Vert B^\star \tilde z + d\Vert_U\Vert B^\star \tilde z\Vert_U\\
\leq & k(\Vert B^\star \tilde z\Vert_U +\Vert d\Vert_U) \Vert B^\star \tilde z\Vert_U\\
\leq & \frac{3}{2}k\Vert B^\star\Vert_{\mathcal{L}(H,U)}^2\Vert \tilde z\Vert_U^2  + \frac{1}{2}k\Vert d\Vert_U^2,
\end{align*}
where \eqref{abstract-difference} has been used for the first inequality, Item 3. of Definition \ref{def-sat} has been used to get the second one, the triangle inequality has been used in the third one, and Cauchy-Schwarz inequality has been used for the last one. Applying Gr\"onwall Lemma, one obtains
\begin{align*}
\Vert \tilde z(t)\Vert_H^2\leq e^{\frac{3}{2}k\Vert B^\star\Vert_{\mathcal{L}(H,U)}^2 t}\Vert \tilde z_0\Vert^2_H + \frac{1}{2}k\int_0^t \Vert d(\tau)\Vert^2_U d\tau
\end{align*}
where $\tilde z_0$ is the initial condition of $\tilde z$. Since the initial condition of $\tilde z$ is null, it follows
\begin{align*}
\Vert \tilde z(t)\Vert_H\leq  \sqrt{\frac{1}{2}k\int_0^t \Vert d(\tau)\Vert^2_U d\tau}
\end{align*}
Using the fact that $\tilde z=z^d-z$ and a triangle inequality yields
\begin{align*}
\Vert z^d\Vert_H\leq \Vert z\Vert_H + \sqrt{\frac{1}{2}k\int_0^t \Vert d(\tau)\Vert^2_U d\tau}
\end{align*}
Now, note that if the origin of \eqref{undisturbed} is globally asymptotically stable, then there exists a class $\mathcal{KL}$ $\alpha$ such that
\begin{equation}
\Vert z\Vert_H\leq \alpha(\Vert z_0\Vert_H,t),\quad \forall t\geq 0.
\end{equation}
Therefore, it implies
\begin{equation}
\Vert z^d\Vert_H\leq \alpha(\Vert z_0\Vert_H,t) + \rho(\Vert d\Vert_{\mathfrak{U}}),
\end{equation}
where $\rho$ is given by: $\rho(s):=\sqrt{\frac{1}{2}k}s.$
\stopSwann
This concludes the proof of Lemma \ref{lemma-estimation}.
\end{proofof}

Thanks to this result, proving ISS for \eqref{sys-dis} reduces to proving global asymptotic stability of the origin of \eqref{nlclosed-loop} for every weak solution to \eqref{nlclosed-loop}. However, when $S\neq U$, the properties given in Definition \ref{def-sat} do not hold when considering weak solutions. Therefore, the following result is needed to link the attractivity for strong and weak solutions. Note that its proof, which relies on a density argument, can be found in \cite{lasiecka2002saturation} and \cite{map2017mcss}.\stopSwann
\begin{lemma}
\label{lemma-density-contraction}
Let $(W_{\sigma}(t))_{t\geq 0}$ be a strongly continuous semigroup of contractions on $H$, a Hilbert space. Let $D(A)$ be dense in $H$. If for all $z_0\in D(A)$, the following holds 
\begin{equation}
\label{attractivity-proof}
\lim_{t\rightarrow +\infty}\Vert W_{\sigma}(t)z_0\Vert_{H}=0,
\end{equation}
then, for all $z_0\in H$,
\begin{equation}
\lim_{t\rightarrow +\infty}\Vert W_{\sigma}(t)z_0\Vert_{H}=0.
\end{equation} 
\end{lemma} 

Thanks to this latter lemma, proving the global asymptotic stability of the origin of \eqref{nlclosed-loop} for every weak solution reduces to proving it for every strong solution to \eqref{nlclosed-loop}. We need another definition of stability before proving the global asymptotic stability.
\begin{definition}
The origin of \eqref{nlclosed-loop} is said to be \emph{semi-globally exponentially stable} if, for any positive value $r$ and any initial condition satisfying $\Vert z_0\Vert_{D(A_\sigma)}\leq r$, there exist two positive values $\mu:=\mu(r)$ and $K:=K(r)$ such that, for every strong solution to \eqref{nlclosed-loop}
\begin{equation}
\Vert W_{\sigma}(t)z_0\Vert_H\leq K e^{-\mu t}\Vert z_0\Vert_H, \quad \forall t\geq 0,
\end{equation} 
where $(W_\sigma(t))_{t\geq 0}$ is the strongly continuous semigroup of contractions generated by the operator $A_\sigma$.
\end{definition}
This definition is inspired by \cite{mcpa2017siam}.  One has also the following result.
\begin{lemma}
\label{semiexp-to-globas}
If the origin of \eqref{nlclosed-loop} is semi-globally exponentially stable, then  it is globally asymptotically stable.
\end{lemma}

\begin{proofof}\textbf{ Lemma \ref{semiexp-to-globas}:}
Consider a strong solution $\tilde{z}$ to \eqref{nlclosed-loop} whose initial condition is such that
\begin{equation}
\Vert \tilde{z}_0\Vert_{D(A)}\leq 1.
\end{equation}
Since the origin of \eqref{nlclosed-loop} is semi-globally exponentially stable, there exist two positive values $\mu_1$ and $K_1$, such that
\begin{equation}
\label{bounded-1}
\Vert W_{\sigma}(t)\tilde{z}_0\Vert_H \leq K_1 e^{-\mu_1 t}\Vert \tilde{z}_0\Vert_H,\quad \forall t\geq 0.
\end{equation}
Now, consider a strong solution $z$ to \eqref{nlclosed-loop} whose initial condition is such that,
\begin{equation}
\Vert z_0\Vert_{D(A)}\leq r,
\end{equation}
where $r$ is a positive value. Hence, since the origin of \eqref{nlclosed-loop} is semi-globally exponentially stable, there exist two positive values $\mu_r$ and $K_r$ such that
\begin{equation}
\label{bounded-r}
\Vert W_{\sigma}(t) z_0\Vert_H\leq K_r e^{-\mu_r t}\Vert z_0\Vert_H,\quad \forall t\geq 0.
\end{equation}
Consequently, setting $T_r=\mu_r^{-1}\ln(rK_r)$, it holds from \eqref{bounded-r}
\begin{equation}
\Vert z_0\Vert_H\leq r\Rightarrow \Vert W_{\sigma}(T_r)z_0\Vert_H \leq e^{-\ln(rK_r)}rK_r=1 
\end{equation}
Therefore, using \eqref{bounded-1} and the fact that $t\mapsto W_{\sigma}(t)$ is a semigroup of contractions, it follows
\begin{align*}
\Vert W_{\sigma}(t)z_0\Vert_H\leq &K_1e^{-\mu_1(t-T_r)}\Vert W_{\sigma}(T_r)z_0\Vert_H,\quad \forall t\geq T_r \\
\leq & K_1e^{\mu_1 T_r}e^{-\mu_1 t}\Vert z_0\Vert_H. 
\end{align*}
This concludes the proof of Lemma \ref{semiexp-to-globas}. 
\end{proofof}

\subsection{Proof of Theorem \ref{thm-abstract-lyap}}

\begin{proofof}\textbf{ Theorem \ref{thm-abstract-lyap}:}
We split the proof of Theorem \ref{thm-abstract-lyap} into two cases. Firstly, we prove Item 1. of Theorem \ref{thm-abstract-lyap} and then Item 2. Indeed, the Lyapunov functions considered in both cases are different. 

\textbf{Case 1: $S=U$.}

Let us consider the following candidate Lyapunov function
\begin{equation}
\begin{split}
V_1(z):=&\langle Pz,z\rangle_H+\frac{2M}{3}\Vert z\Vert_H^3\\
=&V(z) + \frac{2M}{3}\Vert z\Vert_H^3,
\end{split}
\end{equation}
where $P\in \mathcal{L}(H)$ is defined in \eqref{lyapunov-linear-equation} and $M$ is a sufficiently large positive value that will be chosen later. This function, inspired by \cite{liu1996finite}, is clearly positive definite and tends to infinity if the $H$-norm of $z$ does. 

Firstly applying Cauchy-Schwarz inequality, one has, along the strong solutions to (\ref{sys-dis}), 
\begin{align*} 
\frac{d}{dt}V(z)= &\langle Pz,\tilde{A}z\rangle_H+\langle P\tilde{A}z,z\rangle_H \\
&+\langle PB(B^\star z-\sigma(B^\star z),z\rangle_H\\
&+\langle Pz,B(B^\star z-\sigma(B^\star z)\rangle_H\\
& + \langle PB(\sigma(B^\star z)-\sigma(B^\star z +d)),z\rangle_H\\
&+  \langle z,PB(\sigma(B^\star z)-\sigma(B^\star z +d))\rangle_H\\
\leq  & -C\Vert z\Vert^2_H+2\Vert B^\star z\Vert_U\Vert P\Vert_{\mathcal{L}(H)}\Vert B^\star z-\sigma(B^\star z)\Vert_U\\
&+ 2\langle \sigma(B^\star z)-\sigma(B^\star z +d),B^\star Pz\rangle_U,\\
\leq &-C\Vert z\Vert^2_H+2\Vert B^\star z\Vert_U\Vert P\Vert_{\mathcal{L}(H)}\Vert B^\star z-\sigma(B^\star z)\Vert_U\\
&+ 2k\Vert d\Vert_U\Vert B^\star\Vert_{\mathcal{L}(H,U)}\Vert P\Vert_{\mathcal{L}(H)}\Vert z\Vert_H,
\end{align*}
where (\ref{lyapunov-linear-equation}) and Item 3. of Definition \ref{def-sat} have been used to get this inequality. 
Using Item 4. of Definition \ref{def-sat}, Cauchy-Schwarz inequality and the fact that $B^\star$ is bounded in $U$ yields
\begin{align*}
\frac{d}{dt}V(z)\leq & - \left(C-\frac{\Vert B^\star\Vert^2_{\mathcal{L}(H,U)}\Vert P\Vert^2_{\mathcal{L}(H)}}{\varepsilon_1}\right)\Vert z\Vert_H^2 \\
&+ 2\Vert B^\star z\Vert_U\Vert P\Vert_{\mathcal{L}(H)}\Vert B^\star z-\sigma(B^\star z)\Vert_U\\
&+ k^2\varepsilon_1\Vert d\Vert^2_U,\\
\leq & -\left(C-\frac{\Vert B^\star\Vert^2_{\mathcal{L}(H,U)}\Vert P\Vert^2_{\mathcal{L}(H)}}{\varepsilon_1}\right)\Vert z\Vert^2_H\\
&+2\Vert B^\star\Vert_{\mathcal{L}(H,U)}\Vert P\Vert_{\mathcal{L}(H)}\Vert z\Vert_{H}\left(\langle \sigma(B^\star z),B^\star z\rangle_U\right.\\
&\left.+\langle B^\star z,\sigma(B^\star z)\rangle_U\right)\\
&+ k^2\varepsilon_1\Vert d\Vert^2_U\\
\leq & -\left(C-\frac{\Vert B^\star\Vert^2_{\mathcal{L}(H,U)}\Vert P\Vert^2_{\mathcal{L}(H)}}{\varepsilon_1}\right)\Vert z\Vert^2_H\\
&+2\Vert B^\star\Vert_{\mathcal{L}(H,U)}\Vert P\Vert_{\mathcal{L}(H)}\Vert z\Vert_{H}\left(\langle \sigma(B^\star z),B^\star z\rangle_U\right.\\
&\left.+\langle B^\star z,\sigma(B^\star z)\rangle_U\right)+ k^2\varepsilon_1\Vert d\Vert^2_U,
\end{align*}
where $\varepsilon_1$ is a positive value that will be selected later.
Secondly, using the dissipativity of the operator $A_\sigma$, Item 5. of Definition \ref{def-sat} and Young inequality, one has
\begin{align*}
\frac{2M}{3}\frac{d}{dt}\Vert z\Vert_H^3= & M\Vert z\Vert( \langle Az,z\rangle_H + \langle z,Az\rangle_H)\\
&-2M\Vert z\Vert_H\langle B\sigma(B^\star z+d),z\rangle_H\\
\leq &  -2M\Vert z\Vert_H (\langle \sigma(B^\star z),B^\star z\rangle_U\\
&+\langle \sigma(B^\star z)-\sigma(B^\star z +d),B^\star z\rangle_U)\\
\leq &  -2M\Vert z\Vert_H (\langle \sigma(B^\star z),B^\star z\rangle_U)\\
&+2MC_0\Vert z\Vert_H\Vert d\Vert_U\\
\leq & -2M\Vert z\Vert_H (\langle \sigma(B^\star z),B^\star z\rangle_U)\\
&+\frac{2MC_0}{\varepsilon_2}\Vert z\Vert_H^2+2MC_0\varepsilon_2\Vert d\Vert^2_U,
\end{align*}
where $\varepsilon_2$ is a positive value that has to be selected.
Hence, if one chooses $M$, $\varepsilon_1$ and $\varepsilon_2$ as follows
\begin{equation*}
\left\{
\begin{split}
&M\geq 2\Vert B^\star\Vert_{\mathcal{L}(H,U)}\Vert P\Vert_{\mathcal{L}(H)},\\
&\frac{2MC_0}{\varepsilon_2}+\frac{\Vert B^\star\Vert^2_{\mathcal{L}(H,U)}\Vert P\Vert^2_{\mathcal{L}(H)}}{\varepsilon_1}\leq C,
\end{split}
\right.
\end{equation*}
one obtains
\begin{align*}
\label{abs-strict-lyap}
\frac{d}{dt} V_1(z) \leq &- \left(C-\frac{2M}{\varepsilon_2}-\frac{\Vert B^\star\Vert^2_{\mathcal{L}(H,U)}\Vert P\Vert^2_{\mathcal{L}(H)}}{\varepsilon_1}\right)\Vert z\Vert^2_H\\
&+(C_02M\varepsilon_2+k^2\varepsilon_1)\Vert d\Vert^2_U.
\end{align*}
With an appropriate choice of $\varepsilon_1$ and $\varepsilon_2$, it concludes the proof of Item 1. of Theorem \ref{thm-abstract-lyap}.


\textbf{Case 2: $S\neq U$.}

As it has been noticed before, proving ISS reduces to proving the global asymptotic stability of the origin thanks to Lemma \ref{lemma-estimation}. Therefore, in the case where $S\neq U$, we will prove that the origin of \eqref{nlclosed-loop} is semi-globally exponentially stable. Using Lemma \ref{semiexp-to-globas}, it proves that the origin of \eqref{nlclosed-loop} is globally asymptotically stable.

Pick a positive value $r$. 
Consider any strong solution to \eqref{nlclosed-loop} starting from $z_0$ satisfying 
\begin{equation}
\label{initial-condition-bounded}
\Vert z_0\Vert_{D(A)}\leq r.
\end{equation}
Now, focus on the following Lyapunov function
\begin{equation}
\label{lyapunovS}
\begin{split}
V_2(z):=&\langle Pz,z\rangle_H + \tilde{M}r\Vert z\Vert^2_H\\
=&V(z) + \tilde{M}r\Vert z\Vert^2_H
\end{split} 
\end{equation}
where $\tilde{M}$ is a sufficiently large positive value that will be selected later.

Firstly, using Assumption \ref{linearfeedbacklaw}, one has
\begin{align*}
\frac{d}{dt}V(z) = &\langle Pz,\tilde{A}z\rangle_H + \langle P\tilde Az,z\rangle_H \\
&+ \langle PB(B^\star z-\sigma(B^\star z),z\rangle_H\\
&+\langle Pz,B(B^\star z-\sigma(B^\star z)\rangle_H\\
\leq & -C\Vert z\Vert_H^2 + \langle B^\star z-\sigma(B^\star z), B^\star P z\rangle_U\\
&+\langle B^\star P z,B^\star z-\sigma(B^\star z)\rangle_U.
\end{align*}
Since $(U,S)$ is a rigged Hilbert space, one has the following equality
\begin{align*}
&\langle B^\star z-\sigma(B^\star z)), B^\star P z\rangle_U+\langle B^\star P z,B^\star z-\sigma(B^\star z)\rangle_U = \\
&2(B^\star z-\sigma(B^\star z),B^\star P z)_{S^\prime\times S}.
\end{align*}
Hence, applying Cauchy-Schwarz's inequality, one obtains
\begin{align*}
\frac{d}{dt}V(z) \leq -C\Vert z\Vert^2_H + 2\Vert B^\star P z\Vert_S\Vert B^\star z-\sigma(B^\star z)\Vert_{S^\prime}
\end{align*}
Using Item 4. of Definition \ref{def-sat} yields 
\begin{align}\nonumber
\frac{d}{dt}V(z) \leq & - C\Vert z\Vert^2_H \\
&+ 4c_S\Vert Pz\Vert_{D(A)}(\langle B^\star z,\sigma(B^\star z)\rangle_U)\label{eq:12h58},
\end{align}
where, in the second line, it has been used the assumption given in \eqref{continuous-embedding}.

Now, using \eqref{initial-condition-bounded} together with the monotonicity of the function in \eqref{eq:decrease}, along the strong solutions to (\ref{nlclosed-loop}), we have
\begin{equation}
\Vert W_{\sigma}(t)z_0\Vert_{D(A)}\leq \Vert z_0\Vert_{D(A)}. 
\end{equation}
Therefore, noticing that $\Vert Pz\Vert_{D(A)}:=\Vert Pz\Vert_{H}+\Vert PAz\Vert_{H}\leq \Vert P\Vert_{\mathcal{L}(H)}\Vert z\Vert_{D(A)}$, and using (\ref{eq:12h58}) and \eqref{continuous-embedding}, one obtains
\begin{align*}
\frac{d}{dt}V(z) \leq &-C\Vert z\Vert^2_H \\
&+ 4c_Sr \Vert P\Vert_{\mathcal{L(H)}}(\langle z,B\sigma(B^\star z)\rangle_H).
\end{align*}
Moreover, one has, using the dissipativity of the operator $A$\stopSwann
\begin{align*}
\frac{d}{dt}\tilde{M}r\Vert z\Vert^2_H= & \tilde{M}r\langle Az,z\rangle_H + \langle z,Az\rangle_H\\
&-\tilde{M}r\left(\langle z,B\sigma(B^\star z)\rangle_H + \langle B\sigma(B^\star z),z\rangle_H\right)\\
\leq &-4\tilde{M}r\left(\langle z,B\sigma(B^\star z)\rangle_H\right).
\end{align*}
Hence, if one selects $\tilde{M}$ such that
$$
\tilde{M}> 2c_S\Vert P\Vert_{\mathcal{L}(H)}\stopSwann,
$$
it follows
\begin{equation}
\frac{d}{dt}V_2(z)\leq -C\Vert z\Vert^2_H.
\end{equation}
Note that, from \eqref{lyapunovS}, we have, for all $z\in H$,
\begin{equation}
\label{coercivity}
\tilde{M}r \Vert z\Vert^2_H\leq V_2(z) \leq (\Vert P\Vert_{\mathcal{L}(H)}+\tilde{M}r)\Vert z\Vert^2_H
\end{equation}
It yields
\begin{equation}
\begin{split}
\frac{d}{dt}V_2(z)\leq &- \frac{C}{\Vert P\Vert_{\mathcal{L}(H)}+\tilde{M}r}V_2(z)
\end{split}
\end{equation}
Applying Gr\"onwall's inequality, one obtains
\begin{equation}
V_2(W_{\sigma}(t)z_0) \leq e^{-\mu t}V_2(z_0),\quad \forall t\geq 0,
\end{equation}
where $\mu:=\frac{C}{\Vert P\Vert_{\mathcal{L}(H)}+\tilde{M}r}$.

Therefore, using \eqref{coercivity} yields
\begin{equation}
\Vert W_{\sigma}(t)z_0\Vert_H^2 \leq \frac{\Vert P\Vert_{\mathcal{L}(H)}+\tilde{M}r}{\tilde{M}r}e^{-\mu t}\Vert z_0\Vert^2_{H}. 
\end{equation}
Hence, the origin of \eqref{nlclosed-loop} is semi-globally exponentially stable for any strong solution to \eqref{nlclosed-loop}. It yields that the origin is globally asymptotically stable from Lemma \ref{semiexp-to-globas}. Moreover, using Lemma \ref{lemma-density-contraction}, the origin of \eqref{nlclosed-loop} is globally asymptotically for any weak solution to \eqref{nlclosed-loop}. Therefore, using Lemma \ref{lemma-estimation}, this concludes the proof of Item 2. of Theorem \ref{thm-abstract-lyap}.  
\end{proofof}

\section{Example: a Korteweg-de Vries equation}\label{sec:ex}

\subsection{Applying Theorem \ref{thm-abstract-lyap}}

The Korteweg-de Vries equation (for short KdV) describes long waves in water of relatively shallow depth. It has been deeply studied in recent decades  (see e.g., \cite{cerpa2013control} for a nice introduction to the KdV equation in the context of the control). 

The linearized version of the controlled KdV equation can be written as follows
\begin{equation}
\label{KdV}
\left\{
\begin{split}
&z_t + z_x+z_{xxx}=u,\quad (t,x)\in\mathbb{R}_{\geq 0}\times [0,L],\\
&z(t,0)=z(t,L)=0,\quad t\in \mathbb{R}_{\geq 0}\\
&z_x(t,L) = 0,\quad t\in \mathbb{R}_{\geq 0}\\
&z(0,x) = z_0(x),\quad x\in [0,L].
\end{split}
\right.
\end{equation}
Given 
\begin{align*}
&H=U=L^2(0,L),\\
& D(A) :=\lbrace z\in H^3(0,L)\mid z(0)=z(L)=z^\prime(L)=0\rbrace,
\end{align*}
one can describe \eqref{KdV} as an abstract control system with the following operators 
\begin{equation}
A:z\in D(A)\subset H\mapsto -z^\prime-z^{\prime\prime\prime}\in H,\: B=I_H,
\end{equation}
where $I_H$ denotes the identity operator for $H$. A simple computation shows that
$$
\langle Az,z\rangle_H\leq 0,\: \langle A^\star z,z\rangle_H\leq 0.
$$
These inequalities together with the fact that $A$ is a closed operator allow us to apply the L\"umer-Phillips theorem. Therefore, $A$ generates a strongly continuous semigroup of contractions.

If one sets $u=-B^\star z=-z$, then Assumption \ref{linearfeedbacklaw} holds. Indeed, the associated Lyapunov function is 
\begin{equation}
V(z) = \frac{1}{2}\Vert z\vert^2_H.
\end{equation}

Its derivative along the solutions to \eqref{KdV} satisfies
\begin{equation}
\frac{d}{dt} V(z) \leq -V(z). 
\end{equation}
Therefore, the operator $P$ given by \eqref{lyapunov-linear-equation} reduces to the identity.

Now, focus on the case where the control is saturated and where $S\neq U$. Let $S=L^\infty(0,L)$ and $\sigma$ be given by \eqref{satLinfty}. In \cite{mcpa2015kdv_saturating}, the following lemma has been proved
\begin{lemma}[\cite{mcpa2015kdv_saturating}, Lemma 4.]
There exists a positive value $\Delta$ such that,
for all $z\in D(A)$,
\begin{equation}
\Vert z\Vert_{H^1_0(0,L)}\leq \Delta \Vert z\Vert_{D(A)}.
\end{equation}
\end{lemma}
Using the above mentionned result together with the fact that the space $H^1_0(0,L)$ is continuously embedded in $L^\infty(0,L)$, one obtains that
\begin{equation}
\Vert z\Vert_S\leq \Delta \Vert z\Vert_{D(A)}.
\end{equation}
Since $B^\star=I_H$, one has
\begin{equation}
\Vert B^\star z\Vert_{S}\leq \Delta\Vert z\Vert_{D(A)}.
\end{equation}
Therefore, \eqref{continuous-embedding} holds for the linear Korteweg-de Vries equation (\ref{KdV}).

Finally, since all the properties needed to apply Theorem~\ref{thm-abstract-lyap} hold, this proves that the origin of
\begin{equation}
\label{KdV-saturated}
\left\{
\begin{split}
&z_t + z_x+z_{xxx}=-\mathfrak{sat}_{L^\infty(0,L)}(z+d),\: (t,x)\in\mathbb{R}_{\geq 0}\times [0,L],\\
&z(t,0)=z(t,L)=0,\quad t\in \mathbb{R}_{\geq 0}\\
&z_x(t,L) = 0,\quad t\in \mathbb{R}_{\geq 0}\\
&z(0,x) = z_0(x),\quad x\in [0,L],
\end{split}
\right.
\end{equation}
is ISS with respect to $d$.

\subsection{Numerical simulations}
Let us discretize the PDE (\ref{KdV-saturated}) by means of finite difference method (see e.g. \cite{nm_KdV} for an introduction on the numerical scheme of a generalized Korteweg-de Vries equation). The time and the space steps are chosen such that the stability condition of the numerical scheme is satisfied.

We choose $L=2\pi$, $T=9$, $z_0(x)=(1-\cos(x))$ for all $x\in [0,2\pi]$ and $d(t)=0.05\cos(t)$. Let us numerically compute the corresponding solution to (\ref{KdV-saturated}). 

\begin{figure}[h!]
   \begin{minipage}[c]{.50\linewidth}
      \includegraphics[scale=0.6]{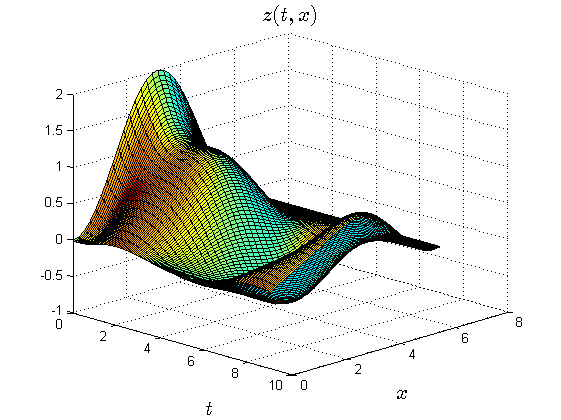}
\end{minipage} \hfill
 \begin{minipage}[c]{.50\linewidth}
      \includegraphics[scale=0.6]{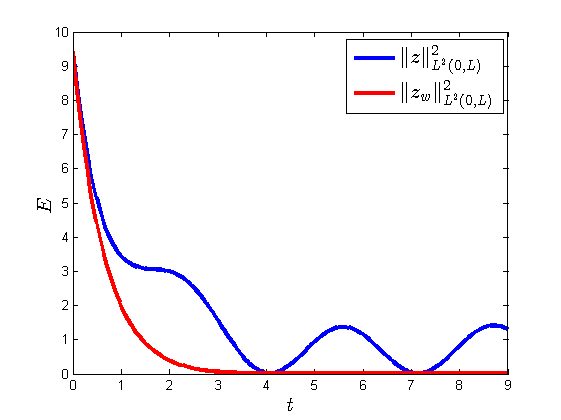}
      \label{figure2}
\end{minipage} \hfill
\caption{Top: time-evolution of the solution to (\ref{KdV-saturated}). Down: $L^2$-norm of the solution to (\ref{KdV-saturated}) (in blue) and of the solution to  (\ref{KdV-saturated}) without saturation and without perturbation (in red).}
\label{figure1}
\end{figure}

Figure \ref{figure1} (top) gives the time-evolution of this solution to \eqref{KdV-saturated}, where the convergence property could be checked. Since the disturbance is bounded, one can check that the solution is also, which is one of the property of ISS. Figure \ref{figure2} illustrates the Lyapunov function $\Vert z\Vert^2_{L^2(0,L)}$ with respect to the time and the solution without perturbation and without saturation. One can see that the $L^2$-norm of the solution to \eqref{KdV-saturated} is bounded.


\section{Conclusion}

\label{sec_conclusion}

In this paper, the analysis of a stabilizing feedback law modified via a saturation function has been tackled using Lyapunov theory for infinite-dimensional systems. The saturation under consideration might render the feedback-law bounded in another space than the one where the origin is stabilized. Assuming a stabilizability property and some regularity for the Lyapunov and the control operators, the system has been proved to be ISS.


For future research lines, it would be interesting to study the case of unbounded control operators as considered for instance in \cite{curtain2006exponential}, \cite{tucsnak2009observation} or \cite{prieur2016wavecone}. The case of nonlinear open-loop systems, as it is considered in \cite{map2017mcss}, could be also challenging. 

\medskip 

\textbf{Acknowledgement:} We would like to thank Professors Marius Tucsnak and Fabian Wirth for fruitful discussions.

\bibliographystyle{IEEEtranS}
\bibliography{bibsm}

\begin{thebibliography}{10}
\providecommand{\url}[1]{#1}
\csname url@samestyle\endcsname
\providecommand{\newblock}{\relax}
\providecommand{\bibinfo}[2]{#2}
\providecommand{\BIBentrySTDinterwordspacing}{\spaceskip=0pt\relax}
\providecommand{\BIBentryALTinterwordstretchfactor}{4}
\providecommand{\BIBentryALTinterwordspacing}{\spaceskip=\fontdimen2\font plus
\BIBentryALTinterwordstretchfactor\fontdimen3\font minus
  \fontdimen4\font\relax}
\providecommand{\BIBforeignlanguage}[2]{{%
\expandafter\ifx\csname l@#1\endcsname\relax
\typeout{** WARNING: IEEEtranS.bst: No hyphenation pattern has been}%
\typeout{** loaded for the language `#1'. Using the pattern for}%
\typeout{** the default language instead.}%
\else
\language=\csname l@#1\endcsname
\fi
#2}}
\providecommand{\BIBdecl}{\relax}
\BIBdecl

\bibitem{azouit2016strong}
R.~Azouit, A.~Chaillet, Y.~Chitour, and L.~Greco, ``Strong i{ISS} for a class
  of systems under saturated feedback,'' \emph{Automatica}, vol.~71, pp.
  272--280, 2016.

\bibitem{cerpa2013control}
E.~Cerpa, ``Control of a {K}orteweg-de {V}ries equation: a tutorial,''
  \emph{Mathematical {C}ontrol and {R}elated {F}ields}, vol. 4(1), pp. 45--99,
  2014.

\bibitem{CHL15}
Y.~Chitour, M.~Harmouche, and S.~Laghrouche, ``${L}_p$-stabilization of
  integrator chains subject to input saturation using {L}yapunov-based
  homogeneous design,'' \emph{SIAM Journal on Control and Optimization},
  vol.~53, no.~4, pp. 2406--2423, 2015.

\bibitem{curtain2006exponential}
R.~Curtain and G.~Weiss, ``Exponential stabilization of well-posed systems by
  colocated feedback,'' \emph{SIAM Journal on Control and Optimization},
  vol.~45, no.~1, pp. 273--297, 2006.

\bibitem{curtain1995introduction}
R.~Curtain and H.~Zwart, \emph{An Introduction to Infinite-Dimensional Systems
  Theory}.\hskip 1em plus 0.5em minus 0.4em\relax Springer-Verlag, New York,
  1995.

\bibitem{Fu}
A.~T. Fuller, ``In-the-large stability of relay and saturating control systems
  with linear controllers,'' \emph{International Journal of Control}, vol.~10,
  no.~4, pp. 457--480, 1969.

\bibitem{karafyllis2016iss}
I.~Karafyllis and M.~Krstic, ``{ISS} with respect to boundary disturbances for
  1-{D} parabolic {PDE}s,'' \emph{IEEE Transactions on Automatic Control},
  vol.~61, no.~12, pp. 3712--3724, 2016.

\bibitem{bible_khalil}
H.~Khalil, \emph{Nonlinear Systems Second Edition}, S.~. Schuster, Ed.\hskip
  1em plus 0.5em minus 0.4em\relax Prentice Hall, Inc., 1996.

\bibitem{lasiecka2002saturation}
I.~Lasiecka and T.~I. Seidman, ``Strong stability of elastic control systems
  with dissipative saturating feedback,'' \emph{Systems \& Control Letters},
  vol.~48, pp. 243--252, 2003.

\bibitem{liu1996finite}
W.~Liu, Y.~Chitour, and E.~Sontag, ``On finite-gain stabilizability of linear
  systems subject to input saturation,'' \emph{SIAM Journal on Control and
  Optimization}, vol.~34, no.~4, pp. 1190--1219, 1996.

\bibitem{map2017mcss}
S.~Marx, V.~Andrieu, and C.~Prieur, ``Cone-bounded feedback laws for
  $m$-dissipative operators on {H}ilbert spaces,'' \emph{Mathematics of
  Control, Signals and Systems}, vol. to appear, 2017.

\bibitem{mcpa2017siam}
S.~Marx, E.~Cerpa, C.~Prieur, and V.~Andrieu, ``Global stabilization of a
  {K}orteweg-de {V}ries equation with a saturating distributed control,''
  \emph{SIAM Journal on Control and Optimization}, vol.~55, no.~3, pp.
  1452--1480, 2017.

\bibitem{mcpa2015kdv_saturating}
------, ``Stabilization of a linear {K}orteweg-de {V}ries with a saturated
  internal control,'' in \emph{Proceedings of the European Control Conference},
  Linz, AU, July 2015, pp. 867--872.

\bibitem{met96}
A.~Megretski, ``${L}_2$ {BIBO} output feedback stabilization with saturated
  control,'' in \emph{Proc. 13th IFAC world congress}, vol. 500, 1996, pp.
  435--440.

\bibitem{MiW16b}
A.~Mironchenko and F.~Wirth, ``Global converse {L}yapunov theorems for
  infinite-dimensional systems,'' in \emph{Proc. of the 10th IFAC Symposium on
  Nonlinear Control Systems (NOLCOS 2016)}, 2016, pp. 909--914.

\bibitem{MiW17b}
------, ``Characterizations of input-to-state stability for
  infinite-dimensional systems,'' \emph{Accepted to IEEE Transactions on
  Automatic Control}, 2017.

\bibitem{nm_KdV}
A.~F. Pazoto, M.~Sep{\'u}lveda, and O.~V. Villagr{\'a}n, ``Uniform
  stabilization of numerical schemes for the critical generalized {K}orteweg-de
  {V}ries equation with damping,'' \emph{Numer. Math.}, vol. 116, no.~2, pp.
  317--356, 2010.

\bibitem{prieur2011iss}
C.~Prieur and F.~Mazenc, ``{ISS} {L}yapunov functions for time-varying
  hyperbolic partial differential equations,'' in \emph{Decision and Control
  and European Control Conference (CDC-ECC), 2011 50th IEEE Conference
  on}.\hskip 1em plus 0.5em minus 0.4em\relax IEEE, 2011, pp. 4915--4920.

\bibitem{prieur2016wavecone}
C.~Prieur, S.~Tarbouriech, and J.~M.~G. da~Silva~Jr, ``Wave equation with
  cone-bounded control laws,'' \emph{IEEE Trans. on Automat. Control}, vol.
  61(11), pp. 3452--3463, 2016.

\bibitem{rao:mag:2001naive}
V.~G. Rao and D.~Bernstein, ``Naive control of the double integrator,''
  \emph{IEEE Control Systems}, vol.~21, no.~5, pp. 86--97, 2001.

\bibitem{Sab0}
A.~Saberi, P.~Hou, and A.~Stoorvogel, ``On simultaneous global external and
  global internal stabilization of critically unstable linear systems with
  saturating actuators,'' \emph{IEEE Transactions on Automatic Control},
  vol.~45, no.~6, pp. 1042--1052, 2000.

\bibitem{seidman2001note}
T.~I. Seidman and H.~Li, ``A note on stabilization with saturating feedback,''
  \emph{Discrete Contin. Dyn. Syst.}, vol.~7, no.~2, pp. 319--328, 2001.

\bibitem{showalter1997monobanach}
R.~Showalter, \emph{Monotone operators in Banach space and nonlinear partial
  differential equations}.\hskip 1em plus 0.5em minus 0.4em\relax Mathematical
  Surveys and Monographs, 1997.

\bibitem{slemrod1989mcss}
M.~Slemrod, ``Feedback stabilization of a linear control system in {H}ilbert
  space with an a priori bounded control,'' \emph{Mathematics of Control,
  Signals and Systems}, vol. 2(3), pp. 847--857, 1989.

\bibitem{sontag1995characterizations}
E.~D. Sontag and Y.~Wang, ``On characterizations of the input-to-state
  stability property,'' \emph{Systems \& Control Letters}, vol.~24, no.~5, pp.
  351--359, 1995.

\bibitem{sontag1989smooth}
E.~Sontag, ``Smooth stabilization implies coprime factorization,'' \emph{IEEE
  transactions on automatic control}, vol.~34, no.~4, pp. 435--443, 1989.

\bibitem{sussmann1991saturation}
H.~Sussmann and Y.~Yang, ``On the stabilizability of multiple integrators by
  means of bounded feedback controls,'' \emph{Technical Report {SYCON}-91-01},
  vol. Rutgers {Center} for {Systems} and {Control}, 1991.

\bibitem{tarbouriech2011book_saturating}
S.~Tarbouriech, G.~Garcia, J.~G. da~Silva~Jr, and I.~Queinnec, \emph{Stability
  and Stabilization of Linear Systems with Saturating Actuators}, S.~Verlag,
  Ed.\hskip 1em plus 0.5em minus 0.4em\relax Springer-Verlag, 2011.

\bibitem{teel1992globalsaturation}
A.~R. Teel, ``Global stabilization and restricted tracking for multiple
  integrators with bounded controls,'' \emph{Systems \& Control Letters},
  vol.~18, pp. 165--171, 1992.

\bibitem{tucsnak2009observation}
M.~Tucsnak and G.~Weiss, \emph{Observation and control for operator
  semigroups}.\hskip 1em plus 0.5em minus 0.4em\relax Springer, 2009.

\end{thebibliography}

\end{document}